# Some remarks on $n_{\max}(k)$ Diophantine $n$-gons

Zurab Aghdgomelashvili

Doctor of Mathematics, Department of Mathematics,
Georgian Technical University
0160, Tbilisi, Georgia

E-mail: z.aghdgomelashvili@gtu.ge

Abstract: In this article one of the fundamental problems of Diophantine (integer) planar geometric figures, the Task $n_{\max}(k)$, is considered.

Task $n_{\max}(k)$. Let us find for each Diophantine $n$-gon $n_0(k) \in N$, for that $n_{\max}(k) \leq n_0(k)$, where $n_{\max}(k)$ – is such maximal $n \in N$, for that the distance between some two vertices of the Diophantine $n$-gon is equal to $k$.

Is determined $n_0(k$ for both convex and concave Diophantine n-gons.

Keywords: convex Diophantine n-gon; concave Diophantine $n$-gon; $n_{\max}(k)$ – is the maximal $n \in N$, at that distance between some two vertices of the Diophantine $n$-gon is equal to $k$.

Introduction.

Let's call as Diophantine a set of points, the distance between arbitrary two points of that is expressed by a natural number.

Let's call as quasi-Diophantine a set of points, the distance between arbitrary two points of that is expressed by a rational number.

One of the most important tasks in the study of Diophantine sets is the following: does exist or not on the Euclidean plane for each number $n \in N$ a having cardinality n Diophantine set, none of three points of that are located on one straight line:

The solution of this problematic task is associated with the name of Wacław Sierpiński. In particular, he proved.

Theorem: for each number $n \in N$ on the Euclidean plane there exists a having cardinality $n$ Quasi-Diophantine set, none of three points of that are located on one straight line.

To solve this task, it was shown that does exist a circumference for that all points of the mentioned Quasi-Diophantine set belongs.

The main idea of constructing the structure of a Quasi-Diophantine set is that initially a dense Quasi-Diophantine set is constructed on a circumference, and if we take a having cardinality $n \in N$ Quasi-Diophantine subset from this set, then with the corresponding homothetic transformation (the coefficient of that, obviously, will be a sufficiently large number), we obtain a Diophantine set of the same cardinality.

This showed that for each $n \in N$ ($n \geq 3$) always exists an $n$-gon, the distance between some two points of that is expressed by a natural number.

This task itself raises serious questions. Among them, one of the most significant is the fundamental task stated by us in the past for Diophantine n-gons, Task ($n; k$). 1.

Task ($n; k$)·1. For each fixed natural number $k$ does exist or not such Diophantine $n$- gon ($n \geq 3$), the distance between some any two vertices of that is equal to $k$ (or is the length of a side or diagonal equal to $k$)? And if it exists, then find all such $n$.

We have shown that does not exist such Diophantine $n$-gon ($n>3$), both convex and concave, the length of arbitrary side or diagonal (distance between some any two vertices) of that is equal to 1, and based on this, only in a Diophantine isosceles triangle the length of the side would be equal to 1. I.e. for $k =1$ $n \in \{3\}$.

We have shown that if $k = 2$, then for convex Diophantine $n$-gons $n \in \{3; 4; 5\}$ (although such a Diophantine pentagon has not yet been found. We assume that such a pentagon does not exist) and for concave $n$-gons $n \in \{3; 4; 5; 6\}$ (although such Diophantine pentagons and hexagons have not yet been found. We assume that such figures do not exist).

We have shown that for Task ($n; 3$)· 1 $3 \leq n \leq 7$. Here we note that for $k=3$ not a single Diophantine pentagon, hexagon or heptagon has been found. We assume that such figures does not exist.

Task ($n; k$) · 1 contains many interesting subtasks. Among them one of most important is the following:

Task $n_{\max}(k)$: Let's find for each Diophantine $n$-gon such $n_0(k) \in N$, for that also $n_{\max}(k) \leq n_0(k)$, where $n_{\max}(k)$ – is the maximal $n \in N$, when the distance between some two vertices of Diophantine tringles $n$-gons is equal to $k$.

Basic part

Task. For solution of $n_{\max}(k)$ let's consider several lemmas.

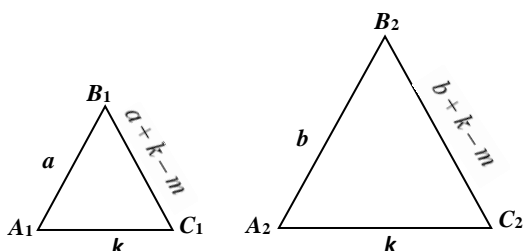

Fig. 1

Lemma 1. Let's consider the stated on Fig. 1 $\Delta A_1B_1C_1$ and $\Delta A_2B_2C_2$.

$a, b, k, m \in N, a < b, \ k > m$

Is necessary to prove that

$$B_1\hat{A}_1C_1 < B_2\hat{A}_2C_2; \quad A_1\hat{C}_1B_1 < A_2\hat{C}_2B_2.$$

Proof: From $\Delta A_1B_1C_1$ and $\Delta A_2B_2C_2$ accordingly of cosine law we have:

$$\cos\hat{C}_1 = \frac{(a+k-m)+k^2-a^2}{2k(a+k-m)}; \quad \cos\hat{C}_2 = \frac{(b+k-m)+k^2-b^2}{2k(b+k-m)};$$

$$\cos\hat{A}_1 = \frac{a^2+k^2-(a+k-m)^2}{2ak}; \quad \cos\hat{A}_2 = \frac{b^2+k^2-(b+k-m)^2}{2bk}.$$

(1)

$(k-m)^2 < k^2 \Leftrightarrow ab + (a+b)(k-m) + (k-m)^2 < k^2 + ab + (k-m)(a+b) \Leftrightarrow$

$\Leftrightarrow (a+k-m)(b+k-m)((a+k-m)-(b+k-m)) < k^2((a+k-m) -$

$- (b+k-m)) + a^2(b+k-m) - b^2(a+k-m) \Leftrightarrow ((a+k-m)^2+k^2-a^2)(b+k-m) >$

$> (a+k-m)((b+k-m)^2+k^2-b^2) \Leftrightarrow \dfrac{(a+k-m)^2+k^2-a^2}{2k(a+k-m)} >$

$> \dfrac{(b+k+m)^2+k^2-b^2}{2k(b+k-m)} \Rightarrow \cos\hat{C}_1 > \cos\hat{C}_2$, and since $\hat{C}_1, \hat{C}_2 \in (0; \pi)$,

thus $\hat{C}_1 < \hat{C}_2$.

Similarly:

$2k > m \Leftrightarrow 2k(b-a) > m^2(b-a) \Leftrightarrow ab(a-b) - k^2(a-b) > b(a+k-m)^2 -$

$- a(b+k-m)^2 \Leftrightarrow b(a^2+k^2-(a+k-m)^2) > a(b^2+k^2-(b+k-m^2)) \Leftrightarrow$

$\Leftrightarrow \dfrac{a^2+k^2-(a+k-m)^2}{2ak} > \dfrac{b^2+k^2-(b+k-m)^2}{2bk} \Rightarrow \cos\hat{A}_1 > \cos\hat{A}_2,$

And as $\hat{A}_1, \hat{A}_2 \in (0; \pi)$, thus $\hat{A}_1 < \hat{A}_2$. Q.E.D.

Consequently we have the following (see Fig. 2).

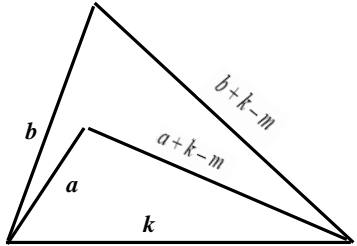

Fig. 2

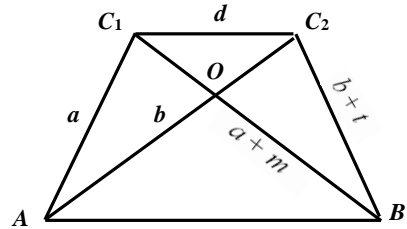

Fig. 3

Lemma 2. Let's show that if for presented on Fig. 3 figure $t, a, b, k, d, m \in N$, then $m > t$.

Proof. From $\Delta AOC_1$ and $\Delta BOC_2$ due the inequality of triangles is easy to obtain:

$$|AC_2| + |C_1B| > |AC_1| + |BC_2| \Rightarrow b + (a + m) > a + (b + t) \Rightarrow m > t. \quad \text{Q.E.D.}$$

Let's consider the Task $n_{max}(k)$ for convex Diophantine $n$-gonsв.

Let's preliminary mention that if $A_1A_2...A_n$ represents the concave $n$-gon, then: the triangle composed of each of its three points does not contain any of the remaining vertices of this Diophantine $n$-gon as an interior point: The value of each of its angles is less than 180°.

Let's consider the case, when the length of side of concave Diophantine $n$-gon is equal to $k$ (Fig.4).

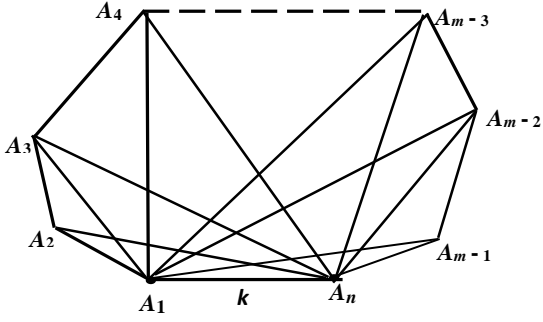

Fig. 4

Based on Lemma 1 and Lemma 2, the cases of maximal $n$ do not exceed the following:

$$|A_2A_n| = |A_1A_2| + k - 1; \ |A_3A_n| = |A_1A_3| + k - 2;$$
$$|A_4A_n| = |A_1A_4| + k - 3;...;|A_mA_n| = |A_1A_m|;...;$$
$$...; |A_1A_{n-1}| = |A_nA_{n-1}| + k - 1.$$

I.e. totally $n_{max}(k) \leq 2(k-1) + 1 + 2 = 2k + 1$.

If $[A_1A_n]$ represents the diagonal, then such $n$-gon would be exist in the second half-plane related to $(A_1A_n)$, and thus:

$$n_{max}(k) \leq 2(2(k-1)+1) + 2 = 4k. \qquad (*)$$

For convex Diophantine $n$-gons prior the consideration of Task $n_{max}(k)$ let's consider the severtal sub-tasks.

Task 1. For the presented on Fig. 5 task let's prove, that $\alpha \geq \beta$.

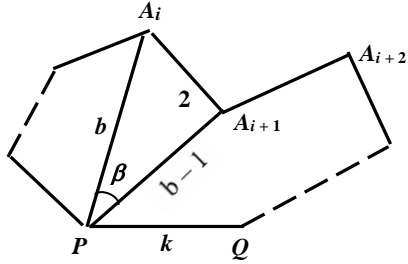 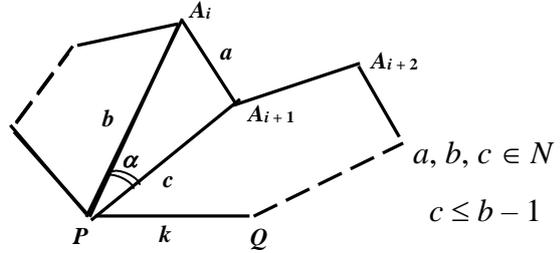

$a, b, c \in N$

$c \leq b - 1$

Fig. 5

Proof: $\cos\beta = \dfrac{2b^2 - 2b - 3}{2b(b-1)}$; $\quad \cos\alpha = \dfrac{b^2 + c^2 - a^2}{2bc}$.

$$c \leq b - 1 \Rightarrow \dfrac{3c}{b-1} \leq 3 \leq (a+1)(a+c-b) \leq (a+b-c)(a+c-b) \Rightarrow$$

$$\Rightarrow \dfrac{b^2 + c^2 - a^2}{2bc} \leq \dfrac{b^2 + (b-1)^2 + 2^2}{2b(b-1)} \Rightarrow \cos\alpha \leq \cos\beta.$$

And as $\alpha, \beta \in (0; \pi)$, thus $\alpha \geq \beta$. Q.E.D.

Task 2. For the presented on Fig. 6 task let's prove, that $\alpha \geq \beta$.

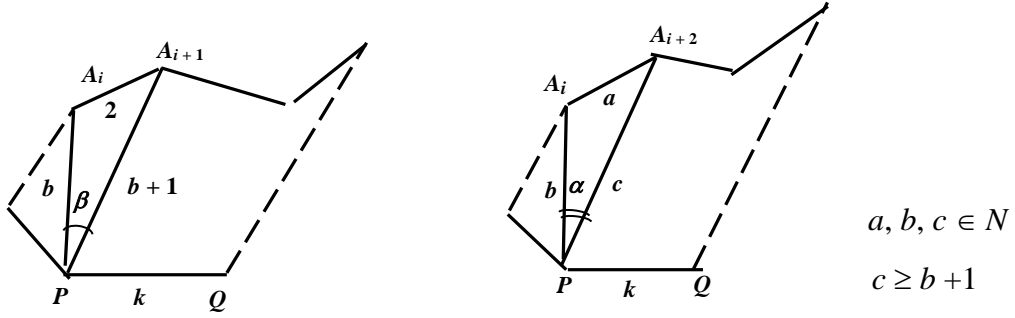

$$a, b, c \in N$$
$$c \geq b + 1$$

Fig. 6

Proof: $\cos\beta = \dfrac{2b(b+1) - 3}{2b(b+1)}$; $\cos\alpha = \dfrac{b^2 + c^2 - a^2}{2bc}$.

Let's say that $c = b + m$ ($m \geq 1$), then the least value of $\alpha$ will be for $a = (m + 1)$. Thus if we show, that $\alpha_{\min} \geq \beta$, then the concept will be proved. We have:

$$\cos\beta = \dfrac{b^2 + (b+1)^2 - 2^2}{2b(b+1)} = 1 - \dfrac{3}{2b(b+1)}; \quad \cos\alpha = \dfrac{b^2 + (b+m)^2 - (m+1)^2}{2b(b+m)} = 1 - \dfrac{2m+1}{2b(b+m)}.$$

$m \geq 2 \Rightarrow (2m-2)b > m - 1 \Rightarrow (2m+1)(b+1) > 3(b+m) \Rightarrow$

$$\Rightarrow \dfrac{2m+1}{2b(b+m)} > \dfrac{3}{2b(b+1)} \Rightarrow 1 - \dfrac{3}{2b(b+1)} > 1 - \dfrac{2m+1}{2b(b+m)} \Rightarrow$$

$$\Rightarrow \dfrac{b^2 + (b+1)^2 - 2^2}{2b(b+1)} > \dfrac{b^2 + (b+m)^2 - (m+1)^2}{2b(b+1)} \Rightarrow \cos\beta > \cos\alpha, \text{ and as } \alpha, \beta \in (0; \pi), \text{ thus } \alpha \geq \beta.$$

Q.E.D.

Task 3. Let's say that from given on the plane $n$ points, none of the three from them are located on the same straight line. Let's arbitrarily select any two of these points, namely $P$ and $Q$ points. Let's show that always does exist an such $n$-gon, each vertex of that belongs to the set of given points, and segment connecting each remaining vertices with $P$ or $Q$ points entirely belongs to this planar $n$-gon.

Proof: To show this, in one half-plane of the $(PQ)$ straight line we draw rays from point $P$ or $Q$ as the origin to the vertices located in this half-plane. By these rays, this half-plane will be divided into non-intersecting parts, each of them is limited by pairs of rays emerging from this origin point. In addition, none of the given points is located in the inner area of each part. We repeat the same procedure would be repeated from the same vertex related to $(PQ)$ in the second half-plane. Let's arbitrarily select any of these rays and denote a point located on it as $A_1$. Then, in clockwise direction, let's designate a point on the adjacent ray as $A_2$, and so on. The last − as $A_n$. Between these vertices will be also the points $P$ and $Q$ (they will be designated with the corresponding numbering). By consequently connecting of these points, we will obtain a polygon $A_1A_2…A_n$ that satisfies the required conditions.

Based on Lemma 1, Lemma 2, Task 1, Task 2 and Task 3, we have that the maximum number of vertices in each half-plane with respect to (PQ) does not exceed $2(k-1)+1 = 2k-1$; therefore, for given Diophantine $n$-gon

$$n_{\max}(k) \leq 2(2k-1)+2 = 4k. \qquad (**)$$

Conclusions

As result of study is determined that: for concave, as well as convex Diophantine $n$-gons we have:

$$n_{\max}(k) \leq 4k. \qquad (***)$$

From this important result follows: for no one $k \in N$ does not found such infinite quantity of points Diophantine sets, fno three points from them are not located on one straight line, and the distance between some any two points is equal to $k$.

## Refferences

.